# Joint Estimation of Multi-phase Traffic Demands at Signalized Intersections Based on Connected Vehicle Trajectories

Chaopeng Tan[a, b], Jiarong Yao[c], Xuegang (Jeff) Ban[d], Keshuang Tang[a, b]*

a.  College of Transportation Engineering, Tongji University, 4800 Cao'an Road, Shanghai 201804, China

b.  Key Laboratory of Road and Traffic Engineering, Ministry of Education, 4800 Cao'an Road, Shanghai 201804, China

c.  School of Electrical and Electronic Engineering, Nanyang Technological University, 50 Nanyang Avenue, Singapore, 639798

d.  The Intelligent Urban Transportation System (iUTS) Lab, Department of Civil and Environmental Engineering, University of Washington, United States

## ABSTRACT

Accurate traffic demand estimation is critical for the dynamic evaluation and optimization of signalized intersections. Existing studies based on connected vehicle (CV) data are designed for a single phase only and have not sufficiently studied the real-time traffic demand estimation for oversaturated traffic conditions. Therefore, this study proposes a cycle-by-cycle multi-phase traffic demand joint estimation method at signalized intersections based on CV data that considers both undersaturated and oversaturated traffic conditions. First, a joint weighted likelihood function of traffic demands for multiple phases is derived given real-time observed CV trajectories, which considers the initial queue and relaxes the first-in-first-out assumption by treating each queued CV as an independent observation. Then, the sample size of the historical CVs is used to derive a joint prior distribution of traffic demands. Ultimately, a joint estimation method based on the maximum a posteriori (i.e., the JO-MAP method) is developed for cycle-based multi-phase traffic demand estimation. The proposed method is evaluated using both simulation and empirical data. Simulation results indicate that the proposed method can produce reliable estimates under different penetration rates, arrival patterns, and traffic demands. The feature of joint estimation makes our method less demanding for the penetration rate of CVs and the consideration of prior distribution can significantly improve the estimation accuracy. Empirical results show that the proposed method achieves accurate cycle-based traffic demand estimation with a MAPE of 12.73%, outperforming the other four methods.

**Keywords**: Traffic demand estimation, connected vehicles, joint estimation, maximum a posteriori, Bayesian theory, oversaturated traffic conditions

# 1 INTRODUCTION

Traffic congestion at urban signalized intersections is usually caused by excessive traffic demand, resulting in the continuous operation of signalized intersections under oversaturated conditions (Koonce and Rodegerdts, 2008; TRB, 2010). Therefore, it is crucial to accurately estimate traffic demands for the dynamic optimization of signal control schemes at intersections. However, most existing studies have focused on signal timing optimization with known traffic demands (Chang and Lin, 2000; Abu-Lebdeh and Benekohal, 2003; Lo and Chow, 2004; Sun et al., 2016), whereas studies on traffic demand acquisition are insufficient, especially for oversaturated traffic conditions (Ma et al., 2017).

Traditionally, traffic demand is measured by fixed-location detectors in traffic control systems deployed at different locations on the roadway, such as stop-bar detectors placed near the stopline, advance detectors placed a few hundred feet upstream from the stopline, and link input detectors placed at the upper reaches of the links (Liu et al., 2009; Hao and Ban, 2015). However, these detectors have limitations in real-world applications, especially in oversaturation scenarios. Stop-bar detectors can only measure the volume of passing vehicles during green phases, that is, the traffic demands for undersaturated traffic conditions. Advance detectors may encounter long queues or spillover problems. In such cases, the detectors are occupied by queued vehicles, and as a result, the detected volume is different from the actual traffic demand. Link input detectors are placed sufficiently upstream from the stopline to measure the input traffic flows of the link while avoiding the influence of long queues. Nevertheless, the exact traffic demands for different lane groups are unknown because of the lane-changing behaviors when vehicles approach the entrance. Although some studies have merged several position detectors to account for the limitations of a single location detector (Ma et al., 2017; Yao and Tang, 2019), their applications are still constrained by the limited detector coverage, malfunctions, and poor data quality of fixed-location detectors (Zhan et al., 2015).

Recent advancements in intelligent mobility and wireless communication have made it possible to access massive amounts of high-resolution vehicle trajectory data from mobile sensors, for example, dedicated short-range communications, smartphones, and other global positioning system-equipped devices, providing great opportunities and unique challenges to enhance traffic monitoring applications. Generally, vehicles equipped with mobile sensors are referred to as connected vehicles (CVs), which can report their instantaneous positions and speeds with a high frequency (i.e., every 1–5 s). Compared with fixed-location sensors that provide traffic information at certain locations, the detection of CVs can almost cover the entire road network and offer a cost-effective data source for monitoring network-wide traffic flows (Ban and Gruteser, 2010).

In the past decade, numerous studies have attempted to estimate performance measures at signalized intersections based on CV data (Guo et al., 2019; Tan et al., 2022). The vast majority of these studies focused on queue length estimation, using either shockwave-based methods (Ban et al., 2011; Cheng et al., 2012; Ramezani and Geroliminis, 2015; Li et al., 2017; Yin et al., 2018), kinematic equation-based methods (Hao et al., 2015), and stochastic methods (Comert and Cetin, 2009, 2011; Hao et al., 2013; Hao et al., 2014; Zhang et al., 2019; Tan et al., 2020; Tan et al., 2021). A few studies have also attempted to estimate the total delay based on CV data. For instance, Ban et al. (2009) estimated delay patterns based on the shockwave theory using travel time data of CVs provided by the virtual trip line technique. Tan et al. (2022) proposed a cumulative flow diagram estimation and prediction method that can simultaneously estimate the average volume, queue length, and total delay of a time-of-day (TOD) period by profiling the cumulative flow diagram of traffic flows based on CV data. Besides, a stream of studies has made efforts to estimate the traffic state (e.g., vehicle counts or traffic density) in the signalized links, using only CV data (Van Phu and Farhi, 2020; Mohammadi and Roncoli, 2021) or fusing CV and fixed-location detector data (Shahrbabaki et al., 2018; Aljamal et al., 2020a; Aljamal et al., 2020b).

Regarding our study object, namely, the estimation of traffic volume (demand), Ban and Gruteser (2010) made the first attempt to estimate cycle-by-cycle traffic volumes at signalized intersections using measured vehicle travel times from mobile sensors. They also aggregated the cycle-by-cycle volumes into longer time durations (from 5 min to 1 h) and showed estimation performance improvement as the aggregation time interval increased. Later studies on this subject can be classified into two categories in terms of the analysis period: TOD-based and cycle-based. Zheng and Liu (2017) are among the first to propose TOD-based traffic volume estimation methods for signalized intersections. Aiming to estimate traffic volumes for short-term periods (e.g., 10, 30, or 60 min) under undersaturated traffic conditions, they formulated the volume estimation as a maximum likelihood estimation problem and solved it using an expectation maximum (EM) algorithm. Zhao et al. (2019) proposed various methods to estimate the period-level maximal queue length and traffic volume by projecting the volume of CVs based on estimated penetration rates but oversaturated traffic conditions were not considered. For cycle-based volume estimation methods, Yao et al. (2020) proposed a traffic volume estimation method for undersaturated cycles that hybridized the shockwave theory and probability distribution. The volume of stopped vehicles was estimated based on the shockwave theory, whereas the estimation of the volume for non-stopped vehicles was modeled as a parameter estimation problem and solved by maximum likelihood estimation. Wang et al. (2019) proposed a framework to combine shockwave analysis with a Bayesian network (BN) for undersaturated traffic state estimation. In their method, a three-layer BN model was applied to construct the relationship between traffic volume, fundamental diagram parameters, arrival times of CVs, and intersection travel time based on shockwave analysis. Furthermore, both cycle-based traffic volumes and queue lengths can be estimated. Later, Wang et al. (2020) extended their method to a pair of adjacent intersections, which achieved queue length and volume estimation for oversaturated traffic conditions. Similar to existing studies using statistical models such as those of Yao et al. (2020) and Zheng and Liu (2017), the BN method also formulated traffic volume estimation into a parameter estimation problem. In contrast to the aforementioned studies, Tang et al. (2020) formulated the problem as a data imputation problem, that is, tensor completion. The traffic volume of each cycle was first divided into two parts, namely, the known part obtained by queued CVs and the unknown part to be estimated. Then, a feature tensor describing cycle-level vehicle arrival during the analysis period was constructed given CV trajectories observed in each cycle, and eventually, the traffic volumes for all cycles were estimated by a tensor decomposition algorithm without any prerequisite assumptions for vehicle arrivals. Nevertheless, oversaturated traffic conditions were not considered in that study.

**TABLE 1 Summary of researches on traffic volume estimation using CV data**

| Research work | Estimated phase | Time resolution | Applicable to oversaturated condition | Real-time or not | Information used* | Methodology |
|---|---|---|---|---|---|---|
| Ban and Gruteser (2010) | Single | Cycle | No | Yes | T | Shockwave analysis |
| Zheng and Liu (2017) | Single | Period | No | No | T, HA | Maximum likelihood estimation |
| Zhao et al. (2019) | Single | Period | No | No | HQP | Queuing position analysis |
| Yao et al. (2019) | Single | Cycle | No | Yes | T | Shockwave analysis & probability distribution |
| Wang et al. (2019) | Single | Cycle | No | No | T | Shockwave analysis & Bayesian network |
| Wang et al. (2020) | Single | Cycle | Yes | No | T | Shockwave analysis & Bayesian network |
| Tang et al. (2020) | Single | Cycle | No | No | T | Tensor completion |
| Tan et al. (2021a) | Single | Period | No | No | T, HA | Total delay calibration |
| This study | Multiple | Cycle | Yes | Yes | T, HA, HQI | Joint maximum a posteriori estimation |

*T – CV trajectory, HA – Historical CV arrivals, HQP – Historical CV queuing positions, HQI – Historical CV quantity information

A summary of relevant research works on traffic volume estimation using CV data is presented in Table 1. Although state-of-the-art methods have achieved cycle-based traffic volume estimation with an acceptable accuracy after years of development, they can still be improved to address the following: 1) All existing studies are designed to estimate the traffic demand of a single phase only. To optimize a signalized intersection, the traffic demand for each phase must be estimated separately. However, owing to the sparsity and randomness of CV sampling in the traffic flow, existing studies for volume estimation of a single phase cannot be easily extended to estimate the volumes of multiple phases because this may cause the traffic demands of some phases to be overestimated, while the others are underestimated. This can lead to a mismatch between the actual demands and the optimized splits obtained by the estimated demands, thereby reducing the efficacy of the intersection operation. 2) Existing studies have not achieved cycle-by-cycle traffic demand estimation for both undersaturated and oversaturated traffic conditions in a real-time manner (estimating the traffic demand of the cycle that just ended). With regard to the category of cycle-based estimation, Yao et al. (2020) and Tang et al. (2020) studied the undersaturated traffic conditions only. Although Wang et al. (2020) studied the oversaturated traffic conditions by analyzing CVs in a pair of adjacent intersections, the method was offline because the three-layer BN model needs to be learned from the CVs of all cycles in a TOD period. 3) Most existing methods only exploit queuing positions, arrival times, or travel time information of CV trajectories, which do not make full use of the quantity information of CVs. However, the number of CVs (i.e., the quantity information) reflects the traffic demands to some extent and has the potential to be used for estimation. The tensor completion method proposed by Tang et al. (2020) adopted the number of CVs as one of the features when constructing the tensor and achieved more accurate estimates than other methods. However, the method is also an offline method that does not consider oversaturated conditions. 4) Existing studies using statistical models, that is, Zheng and Liu (2017) and Yao et al. (2020), commonly used relative queuing and arrival information between consecutive CVs to model vehicle arrivals, which requires the first-in-first-out (FIFO) principle, whereas for multi-lane cases, the FIFO may not be satisfied owing to the variance of vehicle arrivals in different lanes.

To fill the research gap in real-time traffic demand estimation for both undersaturated and oversaturated traffic conditions and address the limitations of existing methods, this study proposes a real-time joint estimation method for multi-phase traffic demands at signalized intersections using CV data. The method considers the prior information provided by the number of CVs based on Bayesian theory, which is applicable to both undersaturated and oversaturated traffic conditions. First, given the queuing and arrival time information of queued CVs during a cycle, a joint weighted likelihood function of traffic demands for multiple phases is derived for the cycle. Because the initial queue is considered during the modeling process, the proposed method also works for oversaturated traffic conditions. In addition, each queued CV is treated as an independent observation of the traffic demand; thus, the FIFO assumption can be relaxed. Furthermore, considering that the sample size of different movements can reflect the relative size of their corresponding traffic demands, the historical observations of CVs for different phases are used to derive the joint prior distribution of traffic demands for multiple phases. Eventually, by maximizing the posterior distribution of traffic demands derived based on the Bayesian theory, the traffic demands of all phases within each cycle can be jointly estimated in a real-time manner.

The major contributions of this study are summarized as follows:

1) This study fills the existing research gap for real-time cycle-based traffic demand estimation for both undersaturated and oversaturated traffic conditions.
2) This study proposes a novel joint modeling method for multi-phase traffic demand estimation based on CV data, which can significantly reduce the requirement for the penetration rate of CVs and meanwhile reduce estimation errors.
3) The integration of the sample size information of CVs into the prior distribution of traffic demands is first proposed in this study, which can improve the estimation accuracy.
4) Two improvements are made when deriving the likelihood function, i.e., treating each queued CV as an independent observation and introducing time-dependent weights to the samples. The former relaxes the FIFO assumption and the latter can improve the reliability of the method.

## 2    PROBLEM DESCRIPTION

As shown in Fig. 1, for the studied phase, the cycles were split by the redstart time. The traffic volume indicates the number of vehicles that traveled through the stopline during the cycle, while the traffic demand indicates the number of vehicles that arrived during the cycle, that is, whose expected arrival time was in the same cycle. Normally, for undersaturated cycles without initial queues, such as cycle $k-2$, the departure vehicles during the green are those vehicles that arrived during the cycle; thus, the traffic demand is the same as the traffic volume, that is, $D_{k-2} = V_{k-2}$. For other cycles that are oversaturated or with initial queues, the traffic demand is different from the traffic volume. For example, cycle $k-1$ is an oversaturated cycle without initial queues; thus, the last few vehicles that arrived in this cycle experienced two stops before traveling through the stopline in the next cycle (i.e., cycle $k$). Therefore, the traffic demand of this cycle is greater than the traffic volume, and it is equal to the traffic volume plus the initial queue of the following cycle, that is, $D_{k-1} = V_{k-1} + q_k^{ini}$. Cycle $k+1$ is an undersaturated cycle with initial queues, and its traffic volume contains the remaining vehicles from the previous cycle, so the traffic demand is less than the traffic volume and equals the traffic volume minus the initial queue, that is, $D_{k+1} = V_{k+1} - q_{k+1}^{ini}$. Cycle $k$ is a combination of cycles $k-1$ and $k+1$; thus, its traffic demand may not be the same as the traffic volume, and we have $D_k = V_k - q_k^{ini} + q_{k+1}^{ini}$.

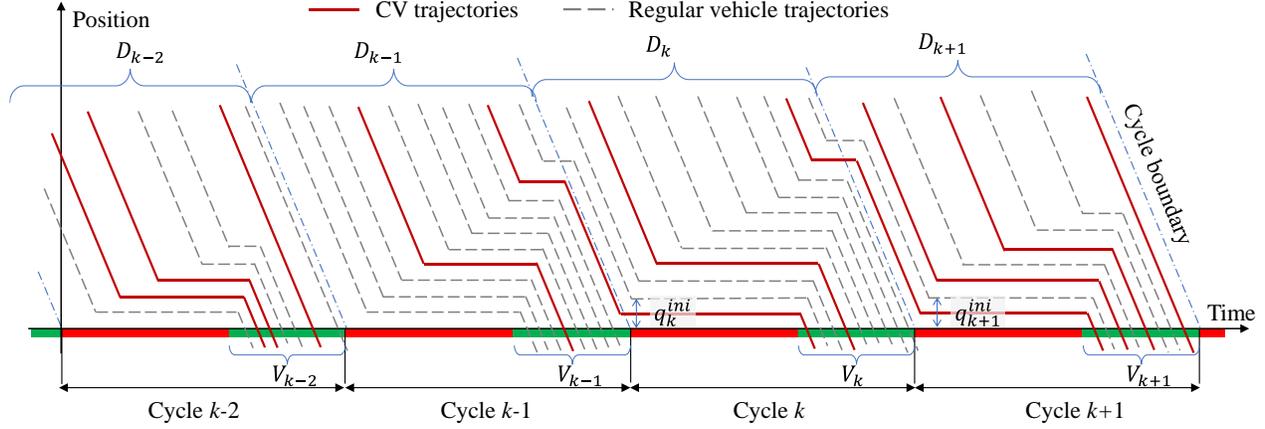

**Fig. 1 Different relationships between traffic demand and traffic volume at signalized intersections**

The traffic demand of a cycle can be expressed in a generalized form and indicated by the average arrival rate, regardless of whether it is oversaturated or if there is an initial queue.

$$D_k = \lambda_k \times C_k \tag{1}$$

where $\lambda_k$ (veh/s) is the average arrival rate of cycle $k$, and $C_k$ (s) is the cycle length of cycle $k$. Notably, in this study, the signal timing of the studied intersection is assumed to be known, and $\lambda_k$ is the only parameter to be estimated.

As mentioned, the number of CVs can reflect the traffic flow demand to some extent. Because the penetration rate is usually unknown in the real world, the number of CV trajectories has little practical use when estimating the traffic demand for a single phase. However, if we broaden the estimation to multiple phases, the number of CV trajectories can provide prior information on the traffic flows by reflecting the relative size of the traffic demands of multiple phases. In light of this, this study exploits the information provided by the number of CV trajectories and proposes a joint estimation method for traffic demands of multiple phases at signalized intersections.

## 3    MODEL FORMULATION

### 3.1    Observations of $\lambda_k$ for a single phase

The cycle a CV belongs to is determined by the cycle where its expected arrival time is located. Given the $j$-th queued CV observed during the $k$-th cycle, we can easily determine its queuing position $n_k^j$ (veh) and expected arrival time during the cycle $t_k^j$ (s), as shown in Eq. (2). The queuing position is defined as the position of a queued vehicle in the queue counting from the stopline, and the expected arrival time refers to the estimated time for the vehicle to arrive at the stopline at a free-flow speed without stopping.

$$n_k^j = \frac{L_k^j}{l_0} \tag{2a}$$

$$t_k^j = T_k^j + \frac{L_k^j}{V_k^j} \tag{2b}$$

where $L_k^j$ is the position where the vehicle joins the queue. Note that, if the vehicle speed is smaller than 5km/h, then the vehicle is identified as queuing status. $l_0$ is the vehicle jam spacing. $T_k^j$ is the time when the vehicle joins the queue. $V_k^j$ is the average speed of the vehicle when approaching the intersection (before joining the queue).

Then, after excluding the initial queue, if it exists, we can prepare a corresponding observation of vehicle arrivals as follows:

$$X_k^j = \left(n_k^j - \hat{q}_k^{ini}, t_k^j\right) \tag{3}$$

where $\hat{q}_k^{ini}$ (veh) is the estimated initial queue of cycle $k$; if not, we have $\hat{q}_k^{ini} = 0$.

According to the flow conservation, the $\tilde{q}_k^{ini}$ is initially calculated, given the estimates of the last cycle:

$$\tilde{q}_k^{ini} = \begin{cases} 0 & if\ \hat{q}_{k-1}^{ini} + \hat{\lambda}_{k-1}C_{k-1} \le s_{k-1}G_{k-1} \\ \hat{q}_{k-1}^{ini} + \hat{\lambda}_{k-1}C_{k-1} - s_{k-1}G_{k-1} & else \end{cases} \tag{4}$$

where $\hat{q}_{k-1}^{ini}$ (veh) is the estimated initial queue of cycle $k-1$, $\hat{\lambda}_{k-1}$ is the estimated average arrival rate of cycle $k-1$, veh/s, $s_{k-1}$ is the saturated departure rate, veh/s, $G_{k-1}$ is the effective green time of cycle $k-1$, s, $q_{k-1}^{ini} + \hat{\lambda}_{k-1}C_{k-1}$ indicates the number of vehicles accumulated during cycle $k-1$, and $s_{k-1}G_{k-1}$ indicates the vehicles dissipated during the effective green time of cycle $k-1$.

$s_{k-1}$ can be calculated based on the departure wave speed fitted by the queue-leaving points of the queued CVs as follows:

$$s_{k-1} = \frac{1}{h_{k-1}^s} = \frac{w_{k-1}^d v_{k-1}^l}{\left(w_{k-1}^d + v_{k-1}^l\right)d^0} \tag{5}$$

where $h_{k-1}^s$ (s) is the saturated time headway of cycle $k-1$, $w_{k-1}^d$ (m/s) is the corresponding departure wave speed, $v_{k-1}^l$ (m/s) is the corresponding average speed value of the queued CVs when they travel through the stopline, and $d^0$ (m) is the empirically calibrated average vehicle jam spacing. Note that when the penetration rate is low, we can calibrate a common $s_{k-1}$ for a TOD period by aggregating all CVs into a common cycle to fit the departure wave.

Referring to the study by Hao et al. (2014), the observed CVs during the last cycle $k-1$ and the current cycle $k$ can help identify whether cycle $k-1$ is oversaturated, and also provide the boundary of the initial queue of cycle $k$. As shown in Fig. 2, CVs are classified into three types according to their queuing process: **Type 1** CVs are vehicles that stop once before passing through the stopline, **Type 2** CVs are vehicles that stop more than once before passing through the stopline, and **Type 3** CVs are vehicles passing through the stopline without stopping.

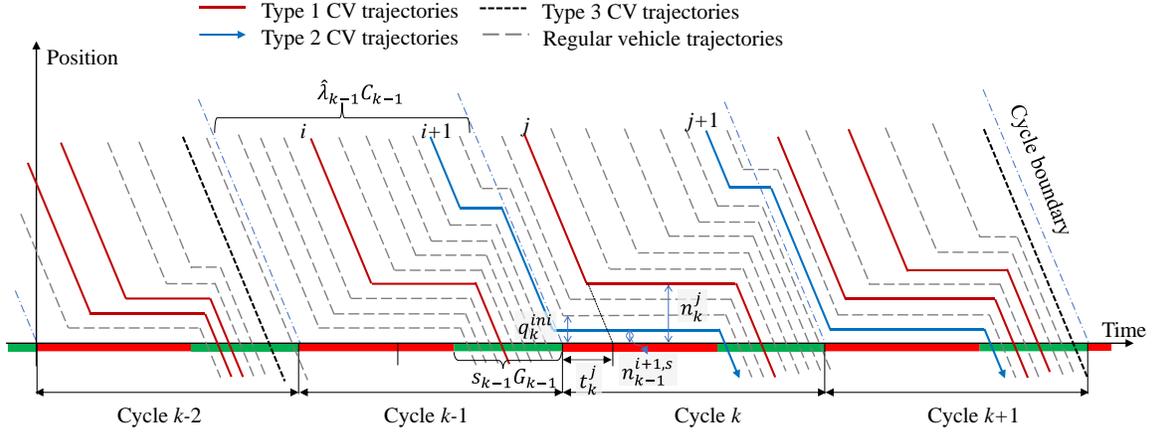

**Fig. 2 Different types of connected vehicles (CVs)**

Different types of CVs can also provide the lower and upper bounds of the initial queue of a cycle. If at least one **Type 2** CV is observed at the last cycle $k-1$, then the initial queue should be greater than the maximum second-time queuing position of **Type 2** CVs, and we have

$$q_k^l = \begin{cases} n_{k-1,2}^{max,s} & \text{if Type 2 CVs are observed at cycle } k-1 \\ 0 & \text{else} \end{cases} \quad (6)$$

where $q_k^l$ (veh) is the lower bound of the initial queue of cycle $k$, and $n_{k-1,2}^{max,s}$ (veh) is the maximum secondary queuing position of **Type 2** CVs from cycle $k-1$.

Concerning the upper bound, if at least one **Type 3** CV is observed in the last cycle, the last cycle is undersaturated, which means that the current cycle $k$ has no initial queue. Then, we have

$$q_k^{u3} = \begin{cases} 0 & \text{if Type 3 CVs are observed at cycle } k-1 \\ \infty & \text{else} \end{cases} \quad (7)$$

where $q_k^{u3}$ (veh) is the upper bound of the initial queue determined by **Type 3** CVs.

**Type 1** CVs can also provide an upper bound of the initial queue if observed.

$$q_k^{u1} = \begin{cases} n_{k,1}^{min} & \text{if Type 1 CVs are observed at cycle } k \\ \infty & \text{else} \end{cases} \quad (8)$$

where $q_k^{u1}$ (veh) is the upper bound of the initial queue determined by **Type 1** CVs, and $n_{k,1}^{min}$ (veh) is the minimum queuing position of **Type 1** CVs in cycle $k$.

Therefore, the final upper bound of the initial queue at cycle $k$ is determined by.

$$q_k^u = \min\left(q_k^{u1}, q_k^{u3}\right) \tag{9}$$

where $q_k^u$ (veh) is the final upper bound of the initial queue in cycle $k$.

Here, we use cycles $k-1$ and $k$ in Fig. 2 as examples to show how we determine the boundary of the initial queues. Because a **Type 3** CV was observed and no **Type 2** CV was observed in cycle $k-2$, and a **Type 1** CV was observed in cycle $k-1$, both the lower bound and the upper bound of the initial queue of this cycle are 0. Concerning cycle $k$, the lower and upper bounds of the initial queue are determined by only **Type 2** CV $i+1$ in cycle $k-1$ and only **Type 1** CV $j$ in cycle $k$; thus, we have $q_k^l = n_{k-1}^{i+1,s}$ and $q_k^u = n_k^j$, where $n_{k-1}^{i+1,s}$ is the secondary queuing position of CV $i+1$ and $n_k^j$ is the queuing position of CV $j$.

In summary, by combining the boundary determined by the observed CVs, the initial queue is eventually estimated as

$$\hat{q}_k^{ini} = \begin{cases} q_k^l & if \ \tilde{q}_k^{ini} < q_k^l \\ \tilde{q}_k^{ini} & if \ q_k^l \leq \tilde{q}_k^{ini} \leq q_k^u \\ q_k^u & if \ \tilde{q}_k^{ini} > q_k^u \end{cases} \tag{10}$$

Each queued CV trajectory can be treated as an independent observation of vehicle arrivals and can be rewritten in the following unified form:

$$X_k^j = \left(n_k^j, t_k^j\right) \tag{11}$$

where $n_k^j$ (veh) indicates the number of vehicles arriving before the $j$-th CV in cycle $k$, $t_k^j$ (s) indicates the corresponding arrival time interval from the start of the cycle to the expected arrival time of the $j$-th CV, and $X_k^j$ indicates that there are $n_k^j$ vehicles arriving during the period $[0, t_k^j]$ in cycle $k$, where the left boundary 0 denotes the relative zero, namely, the start of each cycle.

Considering that similar traffic patterns occur during the same period each day, most signalized intersections operate in a TOD-based mode in the real world. Concurrently, owing to the upstream signals, vehicle arrivals at urban intersections would follow a cyclic pattern in a TOD period, and the cyclic pattern represents the time-varying arrivals within the cycle. To characterize the non-homogeneous arrivals during the cycle, existing studies usually apply a time-dependent Poisson process to model vehicle arrivals (Hao et al., 2013; Hao et al., 2014; Zheng and Liu, 2017; Tan et al., 2021; Tan et al., 2022). The non-homogeneous vehicle arrivals during cycle $k$ can be represented by $\lambda_k \eta(t)$. $\lambda_k$ (veh/s) is the average arrival rate of cycle $k$, and $\eta(t)$ is the time-dependent scaling parameter of the average arrival rate during the cycle. $\eta(t)$ can be obtained by aggregating CV trajectories into a common cycle. For the detailed aggregation process, the reader can refer to Tan et al. (2022). For all cycles during the same TOD period, $\lambda_k$ varies, but $\eta(t)$ is assumed to be the same.

Existing statistical approaches for traffic volume estimation, such as those proposed by Zheng and Liu (2017) and Yao et al. (2020), commonly adopted the relative arrival time and queuing position of consecutive CVs for model construction and assumed that the queue process of vehicles conforms to the FIFO principle. However, this assumption is not always true in the multi-lane case (where the studied movement is served by a lane group). Some vehicles that arrive first may queue at a further queuing position than the later arriving vehicles, which may introduce more errors when modeling vehicle arrivals with consecutive CVs at different lanes in the lane group of the studied phase. To relax this assumption, in this study, all CVs were treated as independent observations of the arrival rate. Therefore, as shown in Fig. 3, given the $j$-th CV observed in the cycle $k$, based on the time-dependent Poisson process, we have $N(0, t_k^j) \sim Poisson(\lambda_k \int_0^{t_k^j} \eta(t) dt)$, and the corresponding probability density function is written as:

$$f(n_k^j, t_k^j | \lambda_k) = \frac{\left[\lambda_k \int_0^{t_k^j} \eta(t) dt\right]^{n_k^j}}{n_k^j!} e^{-\lambda_k \int_0^{t_k^j} \eta(t) dt} \tag{12}$$

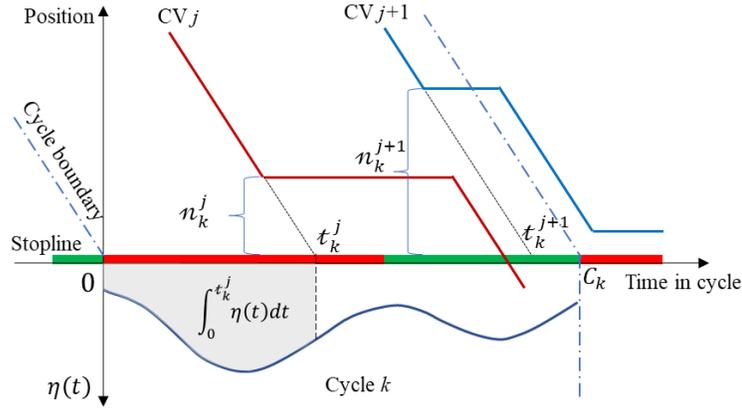

**Fig. 3 Observations of vehicle arrivals based on CVs**

## 3.2 Joint weighted likelihood function

In Section 3.1, we presented how we extracted the observations of vehicle arrivals from CV trajectories for a single phase. In this section, we introduce how to use these observations to establish a joint likelihood function for cycle-based arrival rates for multiple phases.

Without loss of generality, we assume that there are Z ($Z \geq 2$) phases with the same penetration rate at a signalized intersection. The lane-level average arrival rate of phase $z$ in cycle $k$ is $\lambda_{z,k}$, $z = \{1, 2, \dots, Z\}$. Then, we assume variable $\alpha_{z,k}$ as follows:

$$\alpha_{z,k} = \frac{u_z \lambda_{z,k}}{\sum_{\gamma=1}^{Z} u_\gamma \lambda_{\gamma,k}} \quad \forall z \in \{1, 2, \dots, Z\} \tag{13}$$

where $\alpha_{z,k}$ is the proportion of the arrival rate of phase $z$ to the total in cycle $k$, $\sum_{z=1}^{Z} \alpha_{z,k} = 1$, and $u_z$ is the number of lanes of phase $z$.

Let $\sum_{\gamma=1}^{Z} u_\gamma \lambda_{\gamma,k} = \lambda_{0,k}$, that is, the total arrivals of the intersection in cycle $k$, then we have

$$\lambda_{z,k} = \frac{\alpha_{z,k} \lambda_{0,k}}{u_z} \quad \forall z \in \{1, 2, \dots, Z\} \tag{14}$$

We use $x_{z,k}$ (veh) to denote the number of queued CVs of phase $z$ observed in cycle $k$. Based on Equation (12), a weighted likelihood function of $\lambda_{z,k}$ is formulated (Ahmed et al., 2005):

$$L(\lambda_{z,k}) = \prod_{j=1}^{x_{z,k}} f\left(n_{z,k}^j, t_{z,k}^j \middle| \lambda_{z,k}\right)^{w_{z,k}^j} \tag{15}$$

where $w_{z,k}^j$ is the weight of the $j$-th queued CV of phase $z$ observed in cycle $k$.

Compared with the traditional likelihood function which treats each observation equally, the weighted likelihood function considers the reliability of observations and gives them different weights (i.e., $w_{z,k}^j$). In this study, considering the randomness of the distribution of CVs in the traffic flow and the impact of inevitable positioning errors, we believe that the reliability of CV observations is also different, as thus the weighted likelihood function is used.

There are two reasons for the selection of weights in this study. The first reason is that, from a result-oriented perspective, the maximum likelihood estimate is sensitive to the measurement errors of $n_{z,k}^j$ and $t_{z,k}^j$ when $t_{z,k}^j$ is small. For example, for a connected vehicle with a queuing position of 1 veh and an expected arrival time of 2 s, an error of 1 s of the expected arrival time will cause the estimate to become 1 (1/1) or 0.33 (1/3) and an error of 1 veh of the queuing position will cause the estimate to become 1 (2/2) or 0 (0/2). However, if the expected arrival time is 50 s, these errors will have less effect on the estimate. Despite the absolute error of the same size, the difference in the order of magnitude of the denominators will cause a great deviation in the estimates. Since the uploading frequency of CVs ranges from 1s to 5s and the positioning errors of CVs range from 5 to 10 m, the measurement of expected arrival time and queuing position is prone to such random errors. Thus, we would like to give more weight to those CV observations with a greater expected arrival time to reduce the impact of inevitable measurement errors.

The second reason is that, from the observation perspective, a greater expected arrival time $t_{z,k}^j$ covers a larger proportion of the cycle length, indicating a more reliable observation of vehicle arrivals during the cycle (if the expected arrival time $t_{z,k}^j$ is equal to the cycle length, then the observed $n_{z,k}^j$ is exactly the demand of the cycle). To sum up, in our case, we use the following weights of observations:

$$w_{z,k}^j = \int_0^{t_{z,k}^j} \eta_z(t)dt \tag{16}$$

It is worth noting that prior information is also considered in our subsequent estimations. To balance the impact of the prior distribution and likelihood function on the final estimates, the weights need to be normalized as follows:

$$\omega_{z,k}^j = \frac{w_{z,k}^j x_{z,k}}{\sum_{j=1}^{x_{z,k}} w_{z,k}^j} \tag{17}$$

For the arrival rates of the $Z$ phases in cycle $k$, we have the following joint likelihood function:

$$p(\boldsymbol{n_k}, \boldsymbol{t_k}|\lambda_{1,k}, \lambda_{2,k}, \dots \lambda_{Z,k}) = \prod_{z=1}^{Z}\left[\prod_{j=1}^{x_{z,k}} f\left(n_{z,k}^j, t_{z,k}^j|\lambda_{z,k}\right)^{\omega_{z,k}^j}\right]$$ (18)

where $\boldsymbol{n_k}$ is the vector of the queuing position of all CVs in cycle $k$, and $\boldsymbol{t_k}$ is the vector of the expected arrival time of the CVs in cycle $k$.

By substituting $\lambda_{z,k}$ with Equation (14), we have

$$p(\boldsymbol{n_k}, \boldsymbol{t_k}|\boldsymbol{\theta_k}) = \prod_{z=1}^{Z}\left[\prod_{j=1}^{x_{z,k}} f\left(n_{z,k}^j, t_{z,k}^j\left|\frac{\alpha_{z,k}\lambda_{0,k}}{u_z}\right.\right)^{\omega_{z,k}^j}\right]$$ (19)

where $\boldsymbol{\theta_k} = \left[\lambda_{0,k}, \alpha_{1,k}, \alpha_{2,k}, \dots, \alpha_{Z,k}\right]^T$; $\sum_{z=1}^{Z}\alpha_{z,k} = 1$.

### 3.3 Joint prior distribution

Assuming that the variables in $\boldsymbol{\theta_k}$ are independent, the prior distribution $p(\boldsymbol{\theta_k})$ can be written as

$$p(\boldsymbol{\theta_k}) = p(\lambda_{0,k})\prod_{z=1}^{Z} p(\alpha_{z,k})$$ (20)

Because the CV trajectories do not directly provide the prior information on arrival rates, we adopt a uniform distribution as the prior distribution of $\lambda_{0,k}$. Note that, in practical applications, if the historical cycle-based traffic demands of the analysis period are available, we can also use these historical data to generate the prior distribution.

Because the arrival rate for a single flow does not exceed $1/h_s$, the prior distribution of $\lambda_{0,k}$ can be written as

$$p(\lambda_{0,k}) = \begin{cases} \dfrac{h_s}{\sum_{z=1}^{Z} u_z}, 0 < \lambda_0 \leq \dfrac{\sum_{z=1}^{Z} u_z}{h_s} \\ 0, else \end{cases}$$ (21)

As mentioned, the number of CVs can reflect the size of traffic demands to some extent. Thus, for $\alpha_{z,k}$, the prior distribution can be generated by the number of CVs. Here, we use a 2-h TOD period as an example. For phase $z$, we can count the number of CVs every 5 min, as indicated by $N_z^1, N_z^2, \dots, N_z^{24}$. Then, the proportion of the number of CVs of phase $z$ for each 5 min can be obtained as $\frac{N_z^1}{\sum_{z=1}^{Z} N_z^1}, \frac{N_z^2}{\sum_{z=1}^{Z} N_z^2}, \dots, \frac{N_z^{24}}{\sum_{z=1}^{Z} N_z^{24}}$. Then, the prior distribution of $\alpha_{z,k}$ can be generated. In this study, the most commonly used Gaussian distribution was used to fit the prior distribution of $\alpha_{z,k}$. The notation of the prior distribution of $\alpha_{z,k}$ is written as $\alpha_{z,k} \sim \mathcal{N}(\mu_{z,k}, \sigma_{z,k}^2)$ and its probability density function is

$$p(\alpha_{z,k}) = \frac{1}{(2\pi)^{0.5}\sigma_{z,k}} e^{-\frac{(\alpha_{z,k}-\mu_{z,k})^2}{2\sigma_{z,k}^2}}, z = \{1, 2, \dots, Z\}$$ (22)

where $\mu_{z,k}$ and $\sigma_{z,k}^2$ are the mean value and standard derivation of $\frac{N_z^1}{\sum_{z=1}^Z N_z^1}, \frac{N_z^2}{\sum_{z=1}^Z N_z^2}, \cdots, \frac{N_z^{24}}{\sum_{z=1}^Z N_z^{24}}$.

Similar to $\lambda_{0,k}$, if multi-day historical CV data of this TOD period are available, they can also be used to generate the prior distribution. In addition, it should be noted that the prior distribution $p(\boldsymbol{\theta_k})$ is the same for all cycles.

## 3.4  Joint estimation by maximum a posteriori (JO-MAP)

After we deduce the likelihood function and the prior distribution of parameters, the posterior distribution is derived based on the Bayesian theory:

$$p(\boldsymbol{\theta_k}|\boldsymbol{n_k}, \boldsymbol{t_k}) = \frac{p(\boldsymbol{\theta_k})p(\boldsymbol{n_k}, \boldsymbol{t_k}|\boldsymbol{\theta_k})}{p(\boldsymbol{n_k}, \boldsymbol{t_k})} \tag{23}$$

As $p(\boldsymbol{n_k}, \boldsymbol{t_k})$ is independent to $\boldsymbol{\theta_k}$, we have that $p(\boldsymbol{\theta_k}|\boldsymbol{n_k}, \boldsymbol{t_k})$ is proportional to $p(\boldsymbol{\theta_k})p(\boldsymbol{n_k}, \boldsymbol{t_k}|\boldsymbol{\theta_k})$:

$$p(\boldsymbol{\theta_k}|\boldsymbol{n_k}, \boldsymbol{t_k}) \propto p(\boldsymbol{\theta_k})p(\boldsymbol{n_k}, \boldsymbol{t_k}|\boldsymbol{\theta_k}) \tag{24}$$

Thus, based on the maximum a posteriori (MAP) estimation, we have

$$\widehat{\boldsymbol{\theta_k}} = \arg\max p(\boldsymbol{\theta_k}|\boldsymbol{n_k}, \boldsymbol{t_k}) = \arg\max p(\boldsymbol{\theta_k})p(\boldsymbol{n_k}, \boldsymbol{t_k}|\boldsymbol{\theta_k}) \tag{25}$$

The logarithmic posterior distribution is given by

$$\begin{aligned}
\ln p(\boldsymbol{\theta_k})p(\boldsymbol{n_k}, \boldsymbol{t_k}|\boldsymbol{\theta_k}) &= \ln p(\boldsymbol{\theta_k}) + \ln p(\boldsymbol{n_k}, \boldsymbol{t_k}|\boldsymbol{\theta_k}) \\
&= \ln p(\lambda_{0,k}) + \sum_{z=1}^Z \ln p(\alpha_{z,k}) + \sum_{z=1}^Z \sum_{j=1}^{x_{z,k}} \omega_{z,k}^j \ln f\left(n_{z,k}^j, t_{z,k}^j \,\middle|\, \frac{\alpha_{z,k}\lambda_{0,k}}{u_z}\right)
\end{aligned} \tag{26}$$

Note that Equation (25) is subjected to the equality constraint $\sum_{z=1}^Z \alpha_{z,k} = 1$; thus, the Lagrange multiplier method is used to reformulate the original problem into the Lagrangian function:

$$\begin{aligned}
\mathcal{L}(\boldsymbol{\theta_k}, \delta) &= \ln p(\boldsymbol{\theta_k})p(\boldsymbol{n_k}, \boldsymbol{t_k}|\boldsymbol{\theta_k}) - \delta\left(\sum_{z=1}^Z \alpha_{z,k} - 1\right) \\
&= \ln p(\lambda_{0,k}) + \sum_{z=1}^Z \left(-\frac{(\alpha_{z,k} - \mu_{z,k})^2}{2\sigma_{z,k}^2} + \ln\frac{1}{(2\pi)^{0.5}\sigma_{z,k}}\right) \\
&\quad + \sum_{z=1}^Z \sum_{j=1}^{x_{z,k}} \omega_{z,k}^j \left(n_{z,k}^j \ln\alpha_{z,k} + n_{z,k}^j \ln\lambda_{0,k} + n_{z,k}^j \ln\frac{\int_0^{t_{z,k}^j} \eta_z(t)dt}{u_z} - \ln n_{z,k}^j! \right. \\
&\quad \left. - \alpha_{z,k}\lambda_{0,k}\frac{\int_0^{t_{z,k}^j} \eta_z(t)dt}{u_z}\right) - \delta\left(\sum_{z=1}^Z \alpha_{z,k} - 1\right)
\end{aligned} \tag{27}$$

where $\delta$ is the Lagrange multiplier and $\delta \neq 0$.

Because $-\frac{(\alpha_{z,k}-\mu_{z,k})^2}{2\sigma_{z,k}^2}$, $\ln \alpha_{z,k}$, $\ln \lambda_{0,k}$, $-\alpha_{z,k}\lambda_{0,k}$, and $-\delta\alpha_{z,k}$ are concave, the global maximum of $\mathcal{L}(\boldsymbol{\theta_k}, \delta)$ exists.

Aiming to maximize $\mathcal{L}(\boldsymbol{\theta_k}, \delta)$, the MAP estimator of $\boldsymbol{\theta_k}$ can be obtained by solving

$$\nabla_{\boldsymbol{\theta_k}, \delta}\mathcal{L}(\boldsymbol{\theta_k}, \delta) = 0 \tag{28}$$

To summarize, $\boldsymbol{\theta_k}$ and $\delta$ can be estimated by solving $Z + 2$ equations, $Z + 1$ of which are nonlinear:

$$\frac{\partial \mathcal{L}(\boldsymbol{\theta_k}, \delta)}{\partial \lambda_{0,k}} = \sum_{z=1}^{Z} \sum_{j=1}^{x_{z,k}} \omega_{z,k}^j \left[ \frac{n_{z,k}^j}{\lambda_{0,k}} - \frac{\alpha_{z,k} w_{z,k}^j}{u_z} \right] = 0 \tag{29a}$$

$$\frac{\partial \mathcal{L}(\boldsymbol{\theta_k}, \delta)}{\partial \alpha_{z,k}} = -\frac{\alpha_{z,k} - \mu_{z,k}}{\sigma_{z,k}^2} + \sum_{j=1}^{x_{z,k}} \omega_{z,k}^j \left[ \frac{n_{z,k}^j}{\alpha_{z,k}} - \frac{\lambda_{0,k} w_{z,k}^j}{u_z} \right] - \delta = 0, \;\; z = 1,2,\dots,Z \tag{29b}$$

$$\frac{\partial \mathcal{L}(\boldsymbol{\theta_k}, \delta)}{\partial \delta} = \sum_{z=1}^{Z} \alpha_{z,k} - 1 = 0 \tag{29c}$$

By separating variables, these equations can be simplified as

$$\frac{\partial \mathcal{L}(\boldsymbol{\theta_k}, \delta)}{\partial \lambda_{0,k}} = \frac{1}{\lambda_{0,k}} \sum_{z=1}^{Z} \mathcal{N}_{z,k} - \sum_{z=1}^{Z} \mathcal{W}_{z,k} \alpha_{z,k} = 0 \tag{30a}$$

$$\frac{\partial \mathcal{L}(\boldsymbol{\theta_k}, \delta)}{\partial \alpha_{z,k}} = -\frac{\alpha_{z,k}}{\sigma_{z,k}^2} + \frac{\mathcal{N}_{z,k}}{\alpha_{z,k}} - \lambda_{0,k}\mathcal{W}_{z,k} - \delta + \frac{\mu_{z,k}}{\sigma_{z,k}^2} = 0, \;\; z = 1,2,\dots,Z \tag{30b}$$

$$\frac{\partial \mathcal{L}(\boldsymbol{\theta_k}, \delta)}{\partial \delta} = \sum_{z=1}^{Z} \alpha_{z,k} - 1 = 0 \tag{30c}$$

where $\mathcal{N}_{z,k} = \sum_{j=1}^{x_{z,k}} \omega_{z,k}^j n_{z,k}^j$ and $\mathcal{W}_{z,k} = \frac{\sum_{j=1}^{x_{z,k}} \omega_{z,k}^j w_{z,k}^j}{u_z}$. $\mathcal{N}_{z,k}$, $\mathcal{W}_{z,k}$, $\mu_{z,k}$, and $\sigma_{z,k}^2$ are constants determined by the observed CV trajectories and historical CV trajectories, and $u_z$ is a constant indicating the number of lanes. Note that if no queued CVs are observed at a specific phase during cycle $k$, we have $\mathcal{N}_{z,k} = 0$ and $\mathcal{W}_{z,k} = 0$ in Equation (30).

In particular, Equation (30) is a non-linear equation with $Z + 2$ variables, and they may have multiple solutions. In addition, $\lambda_{0,k}$ and $\alpha_{z,k}$ should be nonnegative. Therefore, Equation (30) is further processed to exclude solutions where $\lambda_{0,k}$ and $\alpha_{z,k}$ are negative.

Equation (30b) can be treated as quadratic equations of $\alpha_{z,k}$, so we can rewrite them as

$$-\frac{1}{\sigma_{z,k}^2}\alpha_{z,k}^2 + \left( \frac{\mu_{z,k}}{\sigma_{z,k}^2} - \lambda_{0,k}\mathcal{W}_{z,k} - \delta \right)\alpha_{z,k} + \mathcal{N}_{z,k} = 0, \;\; z = 1,2,\dots,Z \tag{31}$$

By treating $\lambda_{0,k}$ and $\delta$ as constants in Equation (31), we can solve them as

$$\alpha'_{z,k} = \frac{\sigma_{z,k}^2 \left( \frac{\mu_{z,k}}{\sigma_{z,k}^2} - \lambda_{0,k} \mathcal{W}_{z,k} - \delta \right) - \sigma_{z,k}^2 \sqrt{\left( \frac{\mu_{z,k}}{\sigma_{z,k}^2} - \lambda_{0,k} \mathcal{W}_{z,k} - \delta \right)^2 + 4 \frac{\mathcal{N}_{z,k}}{\sigma_{z,k}^2}}}{2} \tag{32a}$$

$$\alpha_{z,k} = \frac{\sigma_{z,k}^2 \left( \frac{\mu_{z,k}}{\sigma_{z,k}^2} - \lambda_{0,k} \mathcal{W}_{z,k} - \delta \right) + \sigma_{z,k}^2 \sqrt{\left( \frac{\mu_{z,k}}{\sigma_{z,k}^2} - \lambda_{0,k} \mathcal{W}_{z,k} - \delta \right)^2 + 4 \frac{\mathcal{N}_{z,k}}{\sigma_{z,k}^2}}}{2} \tag{32b}$$

We can easily find that $\alpha'_{z,k}$ is negative, and $\alpha_{z,k}$ is positive when $\mathcal{N}_{z,k} \geq 0$. Eventually, Equation (30) can be rewritten as

$$\sum_{z=1}^{Z} \mathcal{N}_{z,k} - \lambda_{0,k} \sum_{z=1}^{Z} \mathcal{W}_{z,k} \alpha_{z,k} = 0 \tag{33a}$$

$$\alpha_{z,k} - \frac{\sigma_{z,k}^2 \left( \frac{\mu_{z,k}}{\sigma_{z,k}^2} - \lambda_{0,k} \mathcal{W}_{z,k} - \delta \right) + \sigma_{z,k}^2 \sqrt{\left( \frac{\mu_{z,k}}{\sigma_{z,k}^2} - \lambda_{0,k} \mathcal{W}_{z,k} - \delta \right)^2 + 4 \frac{\mathcal{N}_{z,k}}{\sigma_{z,k}^2}}}{2} = 0 \tag{33b}$$

$$z = 1, 2, \ldots, Z$$

$$\sum_{z=1}^{Z} \alpha_{z,k} - 1 = 0 \tag{33c}$$

By solving Equation (33), we can ensure that the solutions of $\alpha_{z,k}$ ($z = 1, 2, \ldots, Z$) are nonnegative and $\lambda_{0,k}$ is positive.

After we obtain the final estimates, that is, $\widehat{\boldsymbol{\theta}_k}$, the traffic demands during cycle $k$ are calculated as

$$\widehat{D}_{z,k} = \widehat{\lambda}_{0,k} \times \widehat{\alpha}_{z,k} \times \mathcal{C}_k \tag{34}$$

where $\widehat{D}_{z,k}$ (veh) is the estimated traffic demand of phase $z$ during cycle $k$.

## 4 EVALUATION

The performance of the proposed method was evaluated using both simulation data and empirical data. In the simulation evaluation, the proposed method was tested under different traffic demands, arrival patterns, and penetration rates. In the empirical evaluation, two existing state-of-the-art methods of cycle-based traffic volume estimation proposed by Tang et al. (2020) (denoted by the tensor method) and Yao et al. (2020) (denoted by the hybrid method) were also evaluated for comparison.

By maximizing the weighted likelihood function indicated by Equation (10), we can also individually obtain an estimate of the arrival rate of each phase, which does not consider the prior information provided by the number of CVs. We denoted it as the weighted maximum likelihood estimation (WMLE) method.

Given the mean value of $\alpha_{z,k}$ obtained from the number of historical CVs, we can also develop a joint estimation method by maximizing the likelihood function indicated by Equation (19) (i.e., solving Equation (30a)), where $\alpha_{z,k}$ is treated as a constant rather than a variable. This method jointly estimates the arrival rates of multi-phases but does not consider the prior distribution. We denoted it as the joint maximum likelihood estimation. (JO-MLE) method. Hereafter, we will also use the WMLE and the JO-MLE methods for comparison.

The proposed method can jointly estimate the traffic demand of multiple traffic flows using the MAP method, which considers the prior distribution provided by the number of historical CVS. Thus, we use JO-MAP to denote it in the following text.

In simulation evaluation, the prior distribution was generated from the historical CV data, which are generated by randomly sampling the population 10 times with different random seeds. In empirical evaluation, because historical CV data were not available, the prior distribution was generated from the analysis period in the evaluation. Equation (33) is solved by a Python package named "scipy.optimize," which can solve nonlinear equations using numerical methods. The initial values of $\boldsymbol{\theta}_k$ are given by the mean values of the prior distributions. All experiments were conducted on a laptop with a 1.8 GHz 8-core i7-8550U CPU and 16 GB RAM. Generally, it took no more than 3 ms for the estimation of each cycle, which meets the requirements of real-time practical applications.

Three evaluation metrics were used to evaluate the estimation accuracy of the methods, including the mean absolute error (MAE), mean absolute percentage error (MAPE), and success rate (SR), as given by the following equations:

$$MAE = \frac{1}{M} \sum_{m=1}^{M} \left| D_m - \widehat{D}_m \right| \tag{35a}$$

$$MAPE = \frac{1}{M} \sum_{m=1}^{M} \left| \frac{D_m - \widehat{D}_m}{D_m} \right| \tag{35b}$$

$$SR = \frac{S_M}{M} \tag{35c}$$

where $D_m$ is the corresponding ground truth of an estimate, $\widehat{D}_m$ is the estimate, $M$ is the total number of estimations, and $S_M$ is the number of successful estimations of phases and cycles; A successful estimation means that a method can directly obtain an estimate of a certain phase in a certain cycle based on the model, instead of using data repair approaches. The MAE and the MAPE are used to evaluate the accuracy of a method, while the SR is used to evaluate the applicability of a method.

## 4.1 Simulation evaluation

### 4.1.1 Simulation model

As shown in Fig. 4(a), the simulation model was built in VISSIM, based on three consecutive intersections on Jinling Road, Changzhou City, China, where the Jinling–Taihu intersection in the middle was the studied intersection. The intersection layout and signal timing plan of the studied intersection provided by the local traffic management are shown in Fig. 4(b), where the phases are defined according to the NEMA structure. To pursue more realistic scenarios, the truck and turning ratios were collected by the license plate recognition data deployed at the stop lines of three intersections and processed to calibrate the simulation model. The simulation time was 9000 s, in which the first 900 s (warm-up) and the last 900 s (to ensure the integrity of vehicle trajectories) were removed from the analysis. Fig. 4(c) shows the settings of the input demands of eight controlled flows, which correspond to phases in the NEMA structure. As indicated by the demand/capacity (d/c), the phases in group 2 were undersaturated, while those in group 1 were oversaturated on average. Fig. 4(d) shows the ground truth traffic demands of the phases for all cycles, which fluctuate from cycle to cycle. The ground truth traffic demands were obtained from all the trajectories extracted from the VISSIM model. It is worth noting that the right-turn flow on the eastern entrance shares the lane with the straight-through flow; thus, it was added to the straight-through flow. Because of the existence of upstream intersections on Jinling Road, vehicle arrivals on the northern and southern entrances are platoon arrivals, while those on the other entrances are random arrivals. The penetration rates of the CVs were set to 2%–100%. The CVs were randomly sampled from the population of the subject intersection 10 times with different random seeds to ensure the reliability of the results, as thus the actual penetration rates of different phases are slightly different. Normally, the actual average penetration rates of different phases may differ by about 3%, while the real-time penetration rate fluctuates significantly from cycle to cycle.

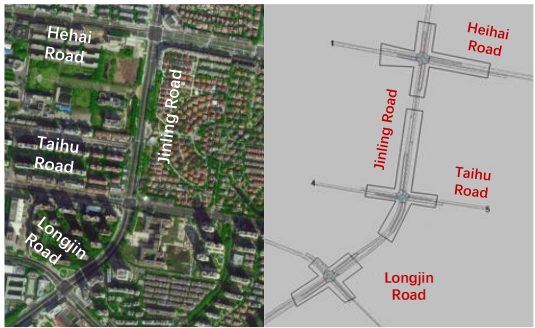

(a) Simulation model of three consecutive intersections

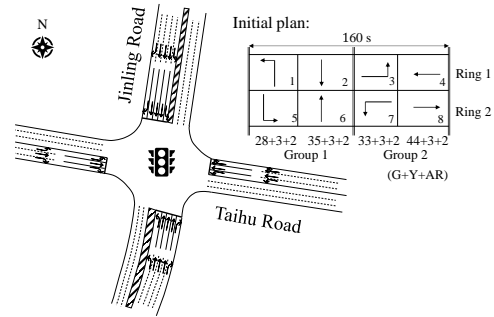

(b) Intersection layout and initial signal timing plan

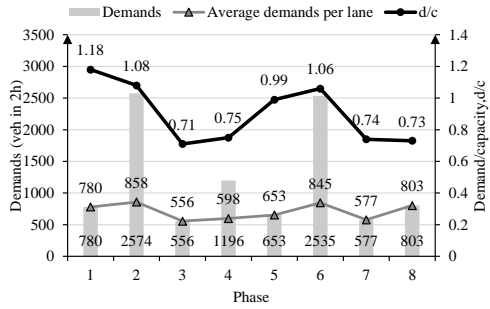

(c) Total demands of phases

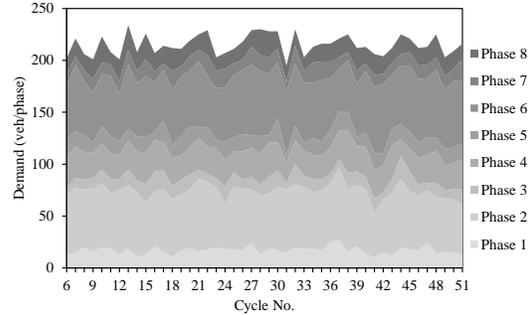

(d) Cycle-based demands of phases

**Fig. 4 Basic information of the simulation model**

### 4.1.2 Different penetration rates

The overall evaluation results of the studied intersections, i.e., the average evaluation results of 8 phases, of WMLE, JO-MLE, and JO-MAP methods at different penetration rates are shown in Fig. 5. With increases in the penetration rate, the MAEs and MAPEs of both the three methods decreased, among which the estimation errors of the proposed JO-MAP methods are always smaller than those of the other two methods. Under low penetration rates (no more than 10%), the JO-MAP method performs slightly better than the JO-MLE method, and significantly better than the WMLE method.

As the penetration rate increases, the error of the WMLE method decreases the fastest, and the MAPE of the WMLE method becomes smaller than that of the JO-MLE method when the penetration rate exceeds 30%. Note that, even with 100% penetration rates, both three methods still have inevitable errors. This is because the methods do not use the penetration rate information, the models do not know that the sample is a full observation and still estimate the arrival rate in a statistical way, which will inevitably bring errors. This is also a common shortage of existing methods that do not rely on penetration rate information.

Fig. 5(c) shows the SRs of the three methods. Because both the JO-MLE method and the JO-MAP method can jointly estimate the traffic demands of multiple phases, their SRs are the same. The SRs of the JO-MAP and the JO-MLE methods are always 100% except for an extremely low penetration rate of 2%, where SR is 93.5%. For the WMLE method, its SR is no more than 90% until the penetration rate increased to 20%, where the SR reaches 94.8%. Note that the SR of the WMLE method indicates the proportion of phases and cycles where queued CVs are observed. When the penetration rate was 2%, queued CVs were observed in only 33.4% of the phases and cycles, while the JO-MAP could still achieve a SR of 95.6%.

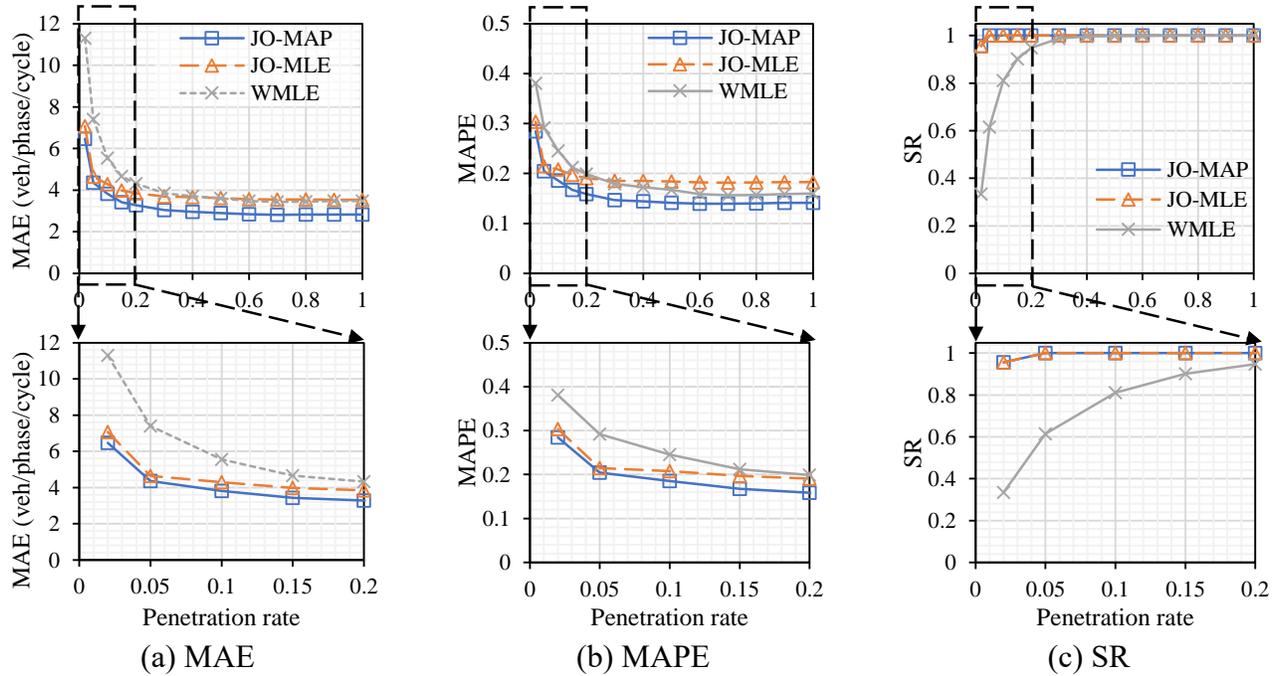

(a) MAE        (b) MAPE        (c) SR

**Fig. 5 Overall performance under different penetration rates**

Fig. 6 presents an example showing the numbers of queued CVs observed in each phase and cycle at 2% and 5% penetration rates. The WMLE method can only be applied to those phases and cycles in which queued CVs are observed (those squares with a number greater than 0). In contrast, the JO-MAP and the JO-MLE methods can jointly estimate the traffic demands of multiple phases. In other words, as long as queued CVs are observed in at least one of the eight phases within a cycle, the JO-MAP and the JO-MLE methods can achieve the demand estimation of all phases of this cycle. Thus, the joint estimation feature of the JO-MAP and the JO-MLE methods can significantly reduce the requirement for the penetration rates of CVs at the intersection.

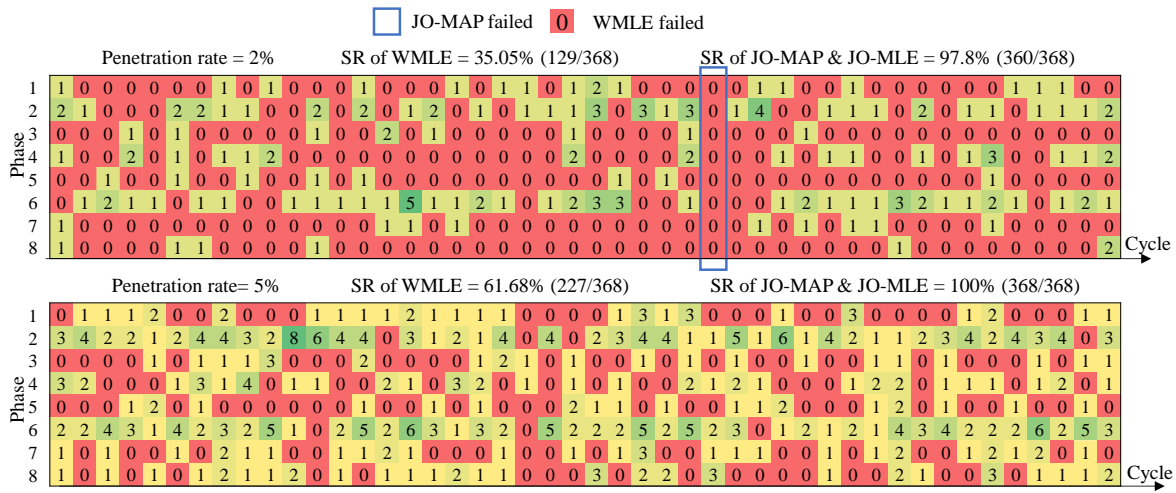

**Fig. 6 Number of queued CVs observed in each phase and cycle (upper figure: 2% penetration rate case; lower figure: 5% penetration rate case)**

Here, we can classify the JO-MAP and the JO-MLE estimates into two parts: estimates of phases and cycles with queued CVs (denoted as with queued CVs) and those without queued CVs (denoted as without queued CVs). Fig. 7 presents the accuracy of the two parts of the JO-MAP and the JO-MLE estimates under different penetration rates. Note that, the corresponding phases and cycles of the results of the JO-MAP and the JO-MLE with queued CVs are the same as WMLE results, i.e., the phases and cycles with successful estimations using the WMLE method.

As the penetration rate increases, the proportion of the phases and cycles with queued CVs increases, and the accuracy of all three estimates, the WMLE, the JO-MAP with queued CVs, and the JO-MLE with queued CVs, increases. However, for those phases and cycles without queued CVs, the accuracy of both the JO-MAP and JO-MLE estimates increases first and then decreases when the penetration rate exceeds a certain value (15% for MAE and 5% for MAPE). This is because when the penetration rate becomes greater, the proportion of the estimates without queued CVs gets smaller. As such estimates mostly occur in those phases and cycles with low traffic demands, so the denominators ($D_m$ and $M$) in Eq.(35) are smaller, resulting in a relatively larger error when calculating the average error of all the estimates without queued CVs. Note that, when the penetration rate is greater than 50%, the SR of the WMLE method reaches 100% (i.e., all phases and cycles can observe queued CVs), therefore there is no estimate for JO-MLE and JO-MAP without queued CVs.

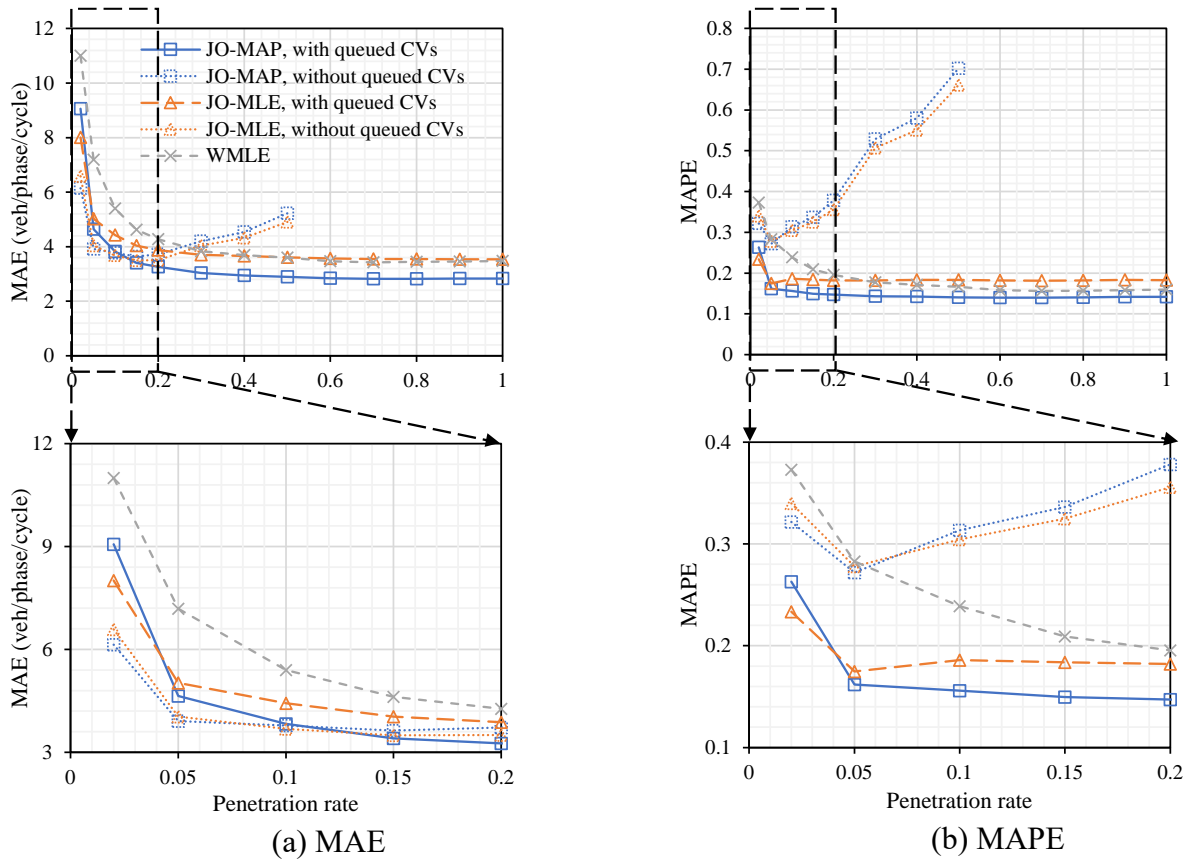

(a) MAE

(b) MAPE

**Fig. 7 Accuracy of JO-MAP estimates with and without queued CVs**

Regarding the estimates of JO-MLE with queued CVs, when the penetration rate is smaller than 30%, the estimates of JO-MLE with queued CVs are more accurate than WMLE by jointly estimating multi-phase traffic demands. However, when the penetration rate exceeds 30%, the estimates of WMLE become more accurate. This is because, in the JO-MLE method, the number of historical CVs is used as a constant. Even a large number of samples is available, the constant still plays a big role to redistribute the results, which brings more errors. In contrast, it can be seen that the estimates of JO-MAP with queued CVs are always more accurate than WMLE. This is because, under a low penetration rate with a small number of samples, the prior distribution can reduce the estimation bias and improve the accuracy; while under a high penetration rate with sufficient samples, the model will rely more on the estimates by samples, and prior distribution only plays a small role to correct the results.

Thus, given the results of JO-MLE and WMLE, we can conclude that by jointly modeling multi-phase traffic demands, the requirement of the penetration rate can be significantly relaxed and the estimation accuracy can be improved under low penetration rates. Given the results of JO-MAP and JO-MLE, we can conclude that by considering the prior distribution generated by historical CVs, the estimation accuracy under all penetration rates can be further improved.

### 4.1.3   Different arrival patterns and d/c

As mentioned before, phases 1, 2, 5, and 6 in group 1 of the NEMA structure are platoon arrivals owing to the existence of the upstream intersection, while the others in group 2 are random arrivals. The overall performance of the two groups under different penetration rates is shown in Fig. 8. Note that because the number of lanes is different for different phases, the evaluation metrics are transformed into the lane average.

When the penetration rate is no more than 40%, the MAEs of group 1 and group 2 are almost the same and the MAPE of group 1 is slightly smaller than group 2. When the penetration rate exceeds 40%, the MAE of group 2 becomes smaller than group 1, while the MAPE of group 2 is slightly greater. Therefore, the difference between the performance of the JO-MAP method in group 1 and group 2 is not significant.

It is worth noting that the phases in Group 1 are all saturated or oversaturated with d/c close to or greater than 1, whereas the phases in Group 2 are undersaturated. In summary, we can conclude that the JO-MAP method is insensitive to arrival patterns and d/c.

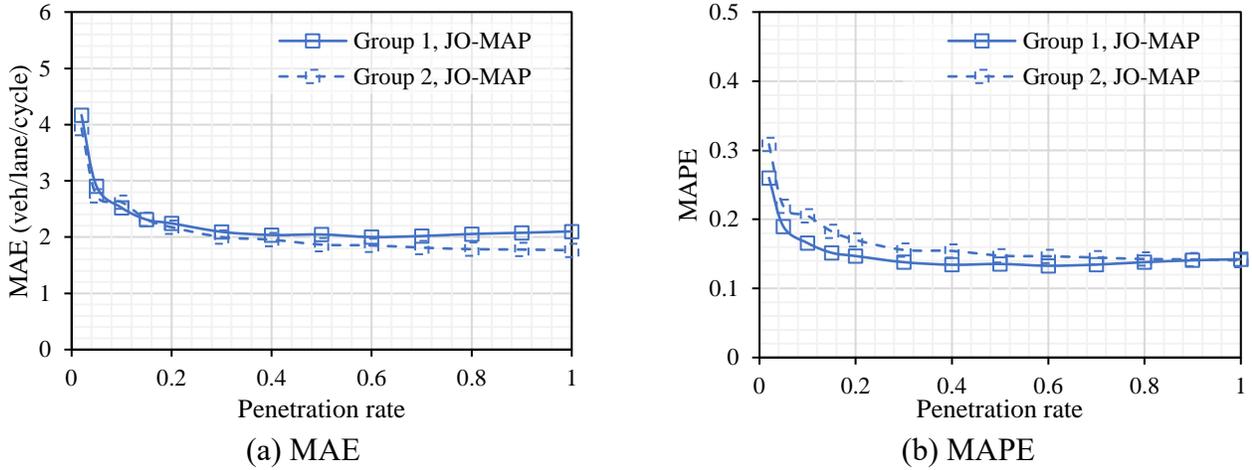

(a) MAE

(b) MAPE

**Fig. 8 Accuracy of JO-MAP estimates with different arrival patterns**

### 4.1.4   *Real-world level penetration rate*

The estimates of the proposed method under a penetration rate of 10%, which is close to the real-world conditions, were used as an example for the analysis. Fig. 9(a) presents a regression plot between the ground truth values and the estimates. From the distribution of the estimates, the WMLE estimates are the most scattered, and some estimates deviate significantly from the ground truth. It is difficult to directly distinguish whether the estimates of the JO-MLE method or the JO-MAP method is more concentrated. . However, the regression results show that the slope of the JO-MAP line is 0.9978, which is the closest to 1 among the three methods, and its R-square value reaches 0.9235, which is also the greatest.

As shown, estimates can be classified into two groups according to ground-truth traffic demands, namely, the low-demand group, whose demands are smaller than 40 veh/phase/cycle, and the high-demand group, whose demands are greater than 40 veh/phase/cycle. The JO-MAP estimates outperformed those of the other two methods with both low and high traffic demands. Regarding the JO-MAP estimates, the MAE of the low-demand group is slightly smaller than that of the high-demand group, whereas the MAPE results are reversed.

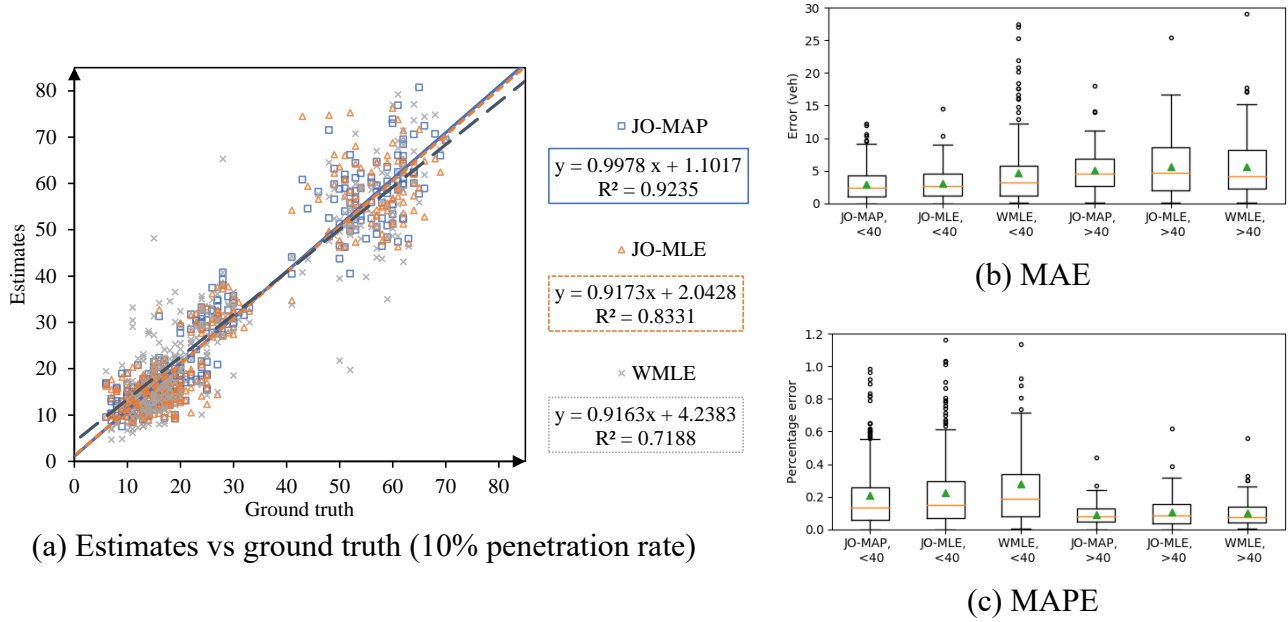

(a) Estimates vs ground truth (10% penetration rate)

(b) MAE

(c) MAPE

**Fig. 9 Accuracy of JO-MAP estimates with a 10% penetration rate**

## 4.2 Empirical evaluation

### 4.2.1 Study site

The Fuzhong Road–Huanggang Road intersection, located in the urban area of Shenzhen, China, was selected as the study site, as shown in Fig. 10(a) and (b). Because only the signal timing and the ground truth demands of the southbound through movement were collected, the four through-lanes in the southbound approach were analyzed in the empirical case. The CV trajectories were provided by Didi Chuxing, and the corresponding ground truth traffic demands for the studied lanes were recorded using a high-resolution video. The uploading interval of CV trajectories was 3 s in most cases, and all data were collected from 10:30 to 14:30 on April 13, 2017.

Because the studied lanes have no initial queues, the traffic demand is the same as the traffic volume in this case. The signal timing plan, ground truth demand, and penetration rate of each cycle are presented in Fig. 10(c). The intersection was operated by the SMOOTH adaptive signal control system in a TOD mode, where four control schemes existed during the study period. In each control scheme, the cycle length and green duration may vary from cycle to cycle with traffic demand. The average penetration rate during the entire analysis period was approximately 8.6%. In total, 88 cycles were used for the evaluation, of which eight cycles had no queued CVs.

Because only the CV trajectories of a single phase were available, we classified the 88 cycles into four groups, each containing 22 cycles, and jointly estimated the traffic demands of the four groups cycle-by-cycle.

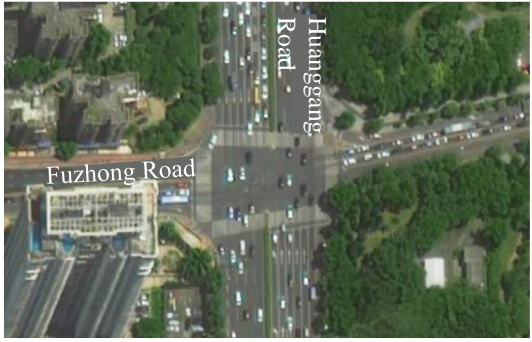

(a) Simulation model of three consecutive intersections

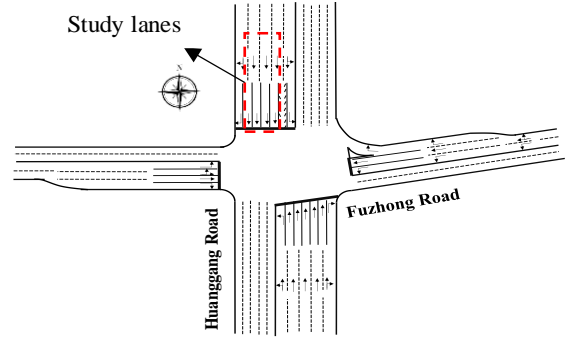

(b) Intersection layout and initial signal timing plan

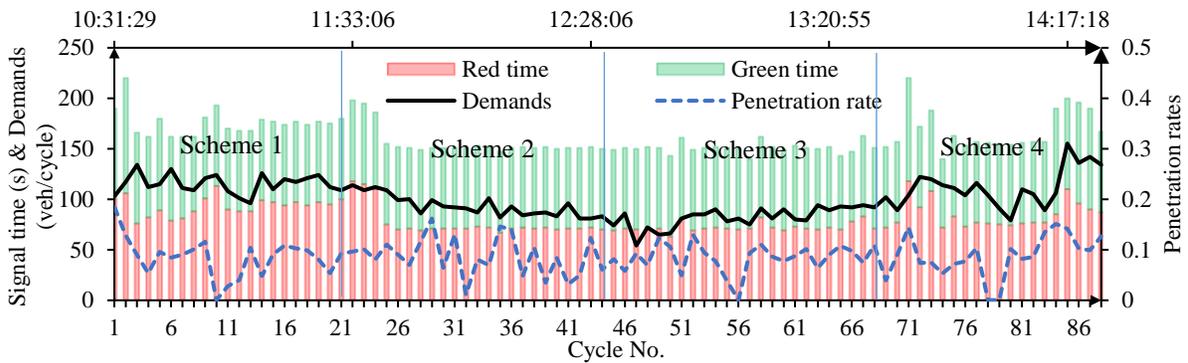

(c) Signal timing plan, ground truth demand, and penetration rate.

**Fig. 10 Basic information of the empirical case**

### 4.2.2 Result analysis

The cycle-based estimates of the tensor method, the hybrid method, the WMLE method, the JO-MLE method, and the proposed JO-MAP method are shown in Fig. 11. As is shown, the overall trends of all methods are consistent with that of the ground truth, whereas the hybrid method may produce large deviations in some cycles. This is because the hybrid method uses two consecutive CVs to model the non-queued volumes, whereas in multi-lane cases, overtaking behaviors near the intersection may result in abnormal observations of traffic arrivals when developing the maximum likelihood function. In contrast, the WMLE, the JO-MLE, and the JO-MAP methods treated CVs as independent observations of traffic arrivals, which relaxes the FIFO assumption and can avoid the impact of overtaking behaviors. The tensor method is a data-driven method that does not require the FIFO assumption either.

By comparing the estimates of the WMLE, the JO-MLE, and the JO-MAP methods, we find that by jointly modeling multi-phase traffic demands, the estimates of the JO-MLE and the JO-MAP methods are more accurate than those of WMLE in most cycles. The regression plots of four methods shown in Fig. 11 also indicate that the JO-MLE and the JO-MAP methods perform better than the other two methods, as the slopes of their regression lines and R-squares are closer to 1. When it comes to the comparison between the JO-MAP estimates and the JO-MLE estimates, we can find that the JO-MAP estimates are more accurate than those of the JO-MLE in 52.3% (46/88) of all 88 cycles. The percentage error of the total demands of 88 cycles of the JO-MAP is 2.9%, which is slightly lower than that of JO-MLE (3.5%).

In addition, the signal timing plan of this case varies from cycle to cycle in different control schemes, and the proposed method does not require the joint estimation objects to be at the same intersection, which means that the proposed method is also applicable to phases at different intersections with different signal control types, provided that their penetration rates are similar.

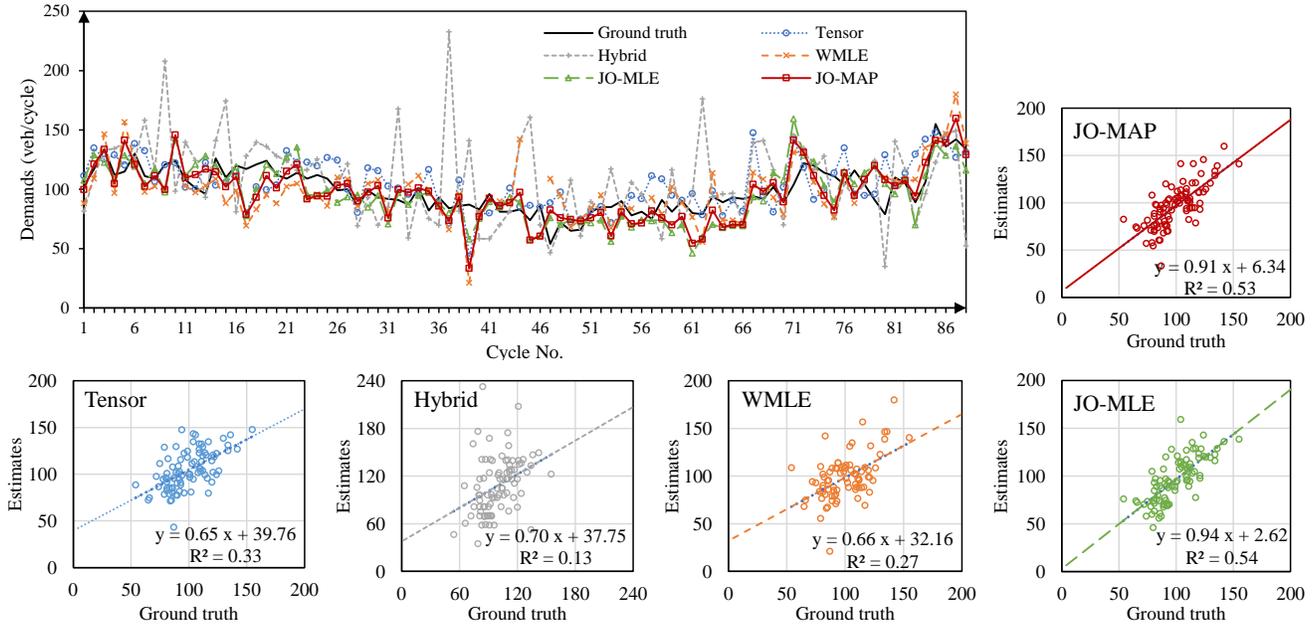

**Fig. 11 Comparisons of cycle-based estimates by different methods**

The overall performances of the five methods are presented in Table 2. Regarding the MAE and MAPE, the proposed JO-MAP method performs best, followed by the JO-MLE, the tensor, and WMLE methods, and the hybrid method performs the worst. Note that, in this case, before we consider the prior information and use the jointly modeling approach, the WMLE method performs worse than the state-of-the-art method, that is, the tensor method. After we use the jointly modeling approach in the JO-MLE method, the accuracy is significantly improved. Furthermore, by integrating the prior distribution based on the Bayesian theory in the JO-MAP method, the accuracy is further improved.

There are eight cycles without queued CVs; thus, the hybrid method and the WMLE method are unavailable for these cycles. Their SRs are 90.9%. The tensor method is a data-driven method that can estimate the volumes of all cycles at once by tensor completion; therefore, its SR is 100%. As for the JO-MLE and the JO-MAP methods, the feature of joint estimation enables it to be applicable even if some cycles lack observations, so its SR is also 100%.

**TABLE 2 Overall performances of the five methods in the empirical case**

| Methods | MAE (veh/cycle) | MAPE | SR |
|---------|-----------------|--------|--------|
| Tensor | 14.35 | 15.27% | 100.0% |
| Hybrid | 24.14 | 25.52% | 90.9% |
| WMLE | 16.75 | 17.42% | 90.9% |
| JO-MLE | 13.18 | 13.80% | 100.0% |
| JO-MAP | 12.73 | 13.37% | 100.0% |

## 5    CONCLUSIONS AND FUTURE STUDIES

To fill the research gap in real-time traffic demand estimation for oversaturated traffic conditions, in this study, we proposed a cycle-by-cycle joint estimation method, named JO-MAP, for traffic demand estimation of multiple phases based on CV data, which works for both undersaturated and oversaturated traffic conditions. In the proposed method, the historical quantities of CVs were used to derive the joint prior distribution of traffic demands of multiple phases, and the real-time observed CVs during the cycle were used to derive the joint weighted likelihood function of traffic demands. Then, based on the Bayesian theory, the cycle-by-cycle traffic demands of multiple phases can be jointly estimated by the maximum posterior distribution. A significant advantage of our method is that we integrate the prior information generated by historical quantities of CVs into our model, which can correct the estimation errors of the likelihood function generated only by observations with a large amount of randomness. Another advantage of our method lies in the joint estimation, which can neutralize the estimation errors of different traffic flows under low penetration rates and significantly reduce the requirement of the minimum penetration rates of CVs, as compared with current methods of estimating traffic flows for a single phase.

The performance of JO-MAP was evaluated using both simulation data and empirical data. Simulation results demonstrate that JO-MAP can provide reliable and accurate estimates under various penetration rates, arrival patterns, and traffic demands. Even with an extremely low penetration rate (2%), where 66.6% of phases in all cycles lack effective observations of CVs, JO-MAP still achieved an SR of 95.7% with an MAE of 6.47 veh/phase/cycle. When the penetration rate is 10%, which is close to the practical application, the MAE of JO-MAP is only 3.82 veh/phase/cycle with an SR of 100%. In the empirical evaluation with an average penetration rate of 8.6%, JO-MAP achieved a MAPE of 13.37% for cycle-based estimates, outperforming two state-of-the-art cycle-based estimation methods with MAPEs of 25.52% and 15.27%, respectively.

The present study also has some limitations that help point out possible directions for future research: 1) Currently, the proposed method assumes that the penetration rates of the studied phases are the same. As this assumption is not strictly satisfied in practical applications, the reliability of the proposed method in scenarios where penetration rates of different phases are significantly different needs to be further evaluated. Future studies are also needed to relax this assumption. 2) Although the phases studied in the evaluations are all at the same intersection, the proposed method does not impose constraints on this. For phases at different intersections with similar penetration rate levels, for example, phases on the mainline of an arterial, the proposed method may also be applicable. Future work may test the proposed method in an arterial or a large-scale road network. 3) The proposed method still requires the signal timing plan to be known, which hinders the application of the proposed method in large-scale road networks. Although some studies have shown that signal timing estimation can be performed for signalized intersections using CV data with relatively high penetration rates (Hao et al., 2012; Li et al., 2017), signal timing estimation with sparse CV data is still challenging. Hence, future work needs to make the proposed method independent of the assumption of known signal timing plans. 4) Future efforts are also needed to extend the proposed method for signalized intersections with shared lanes or storage bays, for which the modeling of traffic flow operation is more complicated and requires further investigation.

## CRediT Authorship Contribution Statement

**Chaopeng Tan**: Conceptualization, Methodology, Validation, Writing - original draft, Writing - review & editing. **Jiarong Yao**: Methodology, Writing - review & editing. **Xuegang (Jeff) Ban**: Conceptualization, Methodology, Writing - review & editing, Supervision. **Keshuang Tang**: Conceptualization, Funding acquisition, Writing - review & editing, Supervision.

## ACKNOWLEDGMENT

This research is supported by the Shanghai Municipal Science and Technology Major Project (2021SHZDZX0100) and the Fundamental Research Funds for the Central Universities. The first author, Chaopeng Tan, gratefully acknowledges financial support from China Scholarship Council (NO: 201906260111) during the visit to the University of Washington.